\documentclass[12pt,a4paper]{article}
\usepackage[dvips]{graphicx}
\usepackage{color}
\usepackage{amssymb}
\usepackage{enumerate}
\usepackage{epsfig}
\usepackage{amsthm}
\usepackage[colorlinks=false,urlcolor=blue,citecolor=blue,linkcolor=blue,bookmarks=true, pdfstartview=FitH, bookmarksopen=true]{hyperref}
\theoremstyle{plain}
\newtheorem{thm}{\noindent Theorem}[section]
\newtheorem{lem}{\noindent Lemma}[section]
\newtheorem{cor}{\noindent Corollary}[section]
\newtheorem{prop}{\noindent Proposition}[section]

\theoremstyle{definition}
\newtheorem{defn}{\noindent Definition}[section]

\theoremstyle{remark}
\newtheorem{rem}{\noindent Remark}[section]

\def\N{\mathbb N}
\def\Z{\mathbb Z}
\def\R{\mathbb R}
\def\C{\mathbb C}
\def\E{\mathbb E}
\def\P{\mathbb P}
\def\1{\textbf{1}}
\def\I{\textrm{1\hspace{-3pt}I}}
\def\rm{\rho_{\max}}
\def\rma{{\mathbb R}_{\max}}
\def\prm{\mathbb{PR}_{\max}^d}
\def\mp{$(\max,+)$}
\def\B{\mathcal B}
\def\M{\mathcal M}

\def\Cr{\mathcal C}

\title{Limit theorems for iterated random topical operators}
\author{Glenn MERLET\\
\href{mailto:glenn.merlet@gmail.com}{glenn.merlet@gmail.com}}

\begin{document}
\maketitle
\begin{abstract}
Let $\left(A(n)\right)_{n\in\N}$ be a sequence of i.i.d. topical (i.e. isotone and additively homogeneous) operators. Let $x(n,x_0)$ be defined by $x(0,x_0)=x_0$ and $x(n,x_0)=A(n)x(n-1,x_0)$.
This can modelize a wide range of systems including, task graphs, train networks, job-shop, timed digital circuits or parallel processing systems.

When $\left(A(n)\right)_{n\in\N}$ has the memory loss property, we use the spectral gap method to prove limit theorems for $x(n,x_0)$. Roughly speaking, we show that $x(n,x_0)$ behaves like a sum of i.i.d. real variables. Precisely, we show that with suitable additional conditions, it satisfies a central limit theorem with rate, a local limit theorem, a renewal theorem and a large deviations principle, and we give an algebraic condition to ensure the positivity of the variance in the CLT. When $A(n)$s are defined by matrices in the \mp semi-ring, we give more effective statements and show that the additional conditions and the positivity of the variance in the CLT are generic.
\end{abstract}

\section*{Introduction}
An operator $A:\R^d\rightarrow\R^d$ is called additively homogeneous if it satisfies $A(x+a\1)=A(x)+a\1$ for all $x\in\R^d$ and $a\in\R$, where $\1$ is the vector $(1,\cdots,1)'$ in $\R^d$. It is called isotone if $x\le y$ implies $A(x)\le A(y)$, where the order is the product order on $\R^d$. It is called topical if it is isotone and homogeneous. The set of topical operators on $\R^d$ will be denoted by $Top_d$.

We recall that the action of matrices with entries in $\rma=\R\cup\{-\infty\}$ on $\rma^d$ is defined by $(Ax)_i=\max_j(A_{ij}+x_j)$. When matrix $A$ has no line of $-\infty$, the restriction of this action to $\R^d$ defines a topical operator, also denoted by $A$. Such operators are called \mp operators and composition of operators corresponds to the product of matrices in the \mp semi-ring.\\

Let $\left(A(n)\right)_{n\in\N}$ be a sequence of random topical operators on $\R^d$. Let $x(n,x_0)$ be defined by 
\begin{equation}\label{defx}
\left\{\begin{array}{lcl}
x(0,x_0)&=&x_0\\
x(n,x_0)&=&A(n)x(n-1,x_0).
\end{array} \right.
\end{equation}

This class of system can modelize a wide range of situations. A review of applications can be found in the last section of~\cite{BM96}. When the $x(n,.)$ are daters, the isotonicity assumption expresses the causality principle, whereas the additive monotonicity expresses the possibility to change the origin of time. (See J.~Gunawardena and M.~Keane~\cite{GunawardenaKeane}, where topical applications have been introduced). Among other examples the \mp case has been applied to modelize queuing networks (J. Mairesse~\cite{Mairesse}, B.~Heidergott~\cite{CaractMpQueuNet}), train networks (B.~Heidergott and R.~De~Vries~\cite{HeidergottDeVriesPubTransNet}, H.~Braker~\cite{braker}) or Job-Shop (G. Cohen and al.~\cite{cohen85a}). It also computes the daters of some task resources models (S.~Gaubert and J.~Mairesse~\cite{gaumair95}) and  timed Petri Nets including Events graphs (F.~Baccelli~\cite{Baccelli}) and 1-bounded Petri Nets (S.~Gaubert and J.~Mairesse~\cite{GaubertMairesseIEEE}). The role of the max operation is synchronizing different events. For devlopements on the max-plus modelizing power, see F.~Baccelli and al.~\cite{BCOQ} or B.~Heidergott, G.~J.~Olsder, and J.~van~der~Woude~\cite{MpAtWork}.

We are interested in the asymptotic behavior of $x(n,.)$. It follows from theorem~\ref{thVincent} that $\frac{1}{n}\max_ix_i(n,X_0)$ converges to a limit $\gamma$.

In many cases, if the system is closed, then every coordinate $x_i(n,X_0)$ also converges to $\gamma$. The value $\gamma$, which is often called cycle time, is the inverse of the throughput (resp. output) of the modelized network (resp. production system), therefore there has been many attempt to estimate it. (J.E.~Cohen~\cite{Cohen}, B.~Gaujal and A.~Jean-Marie~\cite{ComputIssuesSRS}, J.~Resing and al.~\cite{RVH}) Even when the $A(n)$'s are i.i.d. and take only finitely many values, approximating $\gamma$ is NP-hard (V.~Blondel and al.~\cite{LyapExpNP}). D.~Hong and its coauthors have obtained (\cite{BacHong1},\cite{BacHong2} ,\cite{GaubertHong} ) analyticity of $\gamma$ as a function of the law of $A(1)$. In this paper, we prove another type of stability, under the same assumptions.

We show that under suitable additional conditions, $x(n,.)$ satisfies a central limit theorem, a local limit theorem, a renewal theorem and a large deviations principle. When the $A(n)$ are \mp operators we give more explicit results. Those results justify the approximation of $\gamma$ by $\frac{1}{n}x_i(n,X_0)$ and lead the path to confidence intervals.\\

Products of random matrices in the usual sense have been intensively investigated. Let us cite H. Furstenberg~\cite{Furst63}, Y. Guivarc'h and A. Raugi~\cite{GR85} or I. Ya. Gol{$'$}dshe\u{\i}d and G. A. Margulis~\cite{GM2}. The interested reader can find a presentation of this theory in the book by Ph.~Bougerol and J.~Lacroix~\cite{BL}.
We investigate analogous problems to those studied by \'E. Le Page~\cite{LePage}, but for matrices in the \mp semi-ring and more generally for iterated topical operators.

This article is divided into three parts. 
First we present the model of iterated topical operators, including a short review of known limit theorems and a sketch of the proof of our results. Second we state our theorems and comment on them. Finally we prove them.

\section{Iterated topical operators}\label{IFS}
\subsection{Memory loss property}
Dealing with homogeneous operators it is natural to introduce the quotient space of $\R^d$ by the equivalence relation $\sim$ defined by $x\sim y$ if $x-y$ is proportional to $\1$. This space will be called projective space and denoted by $\prm$. Moreover $\overline{x}$ will be the equivalence class of $x$.

The application $\overline{x}\mapsto (x_i-x_j)_{i<j}$ embeds $\prm$ onto a subspace of $\R^{\frac{d(d-1)}{2}}$ with dimension $d-1$. The infinity norm of $\R^{\frac{d(d-1)}{2}}$ therefore induces a distance on $\prm$ which will be denoted by $\delta$. A direct computation shows that $\delta(\overline{x},\overline{y})=\max_i(x_i-y_i)+\max_i(y_i-x_i)$. By a slight abuse, we will also write $\delta(x,y)$ for $\delta(\overline{x},\overline{y})$. The projective norm of $x$ will be $|x|_\mathcal{P}=\delta(x,0)$.

Let us recall two well known facts about topical operators. First a topical operator is non-expanding with respect to the infinity norm. Second the operator it defines from $\prm$ to itself is non-expanding for $\delta$.

The key property for our proofs is the following:
\begin{defn}[MLP]\
\begin{enumerate}
\item A topical operator $A$ is said to have rank~1, if it defines a constant operator on $\prm$ : $\overline{Ax}$ does not depend on $x\in\R^d$.
\item The sequence $\left(A(n)\right)_{n\in\N}$ of $Top_d$-valued random variables is said to have the memory loss (MLP) property  if there exists an $N$ such that $A(N)\cdots A(1)$ has rank~1 with positive probability.
\end{enumerate} 
\end{defn}
This notion has been introduced by J. Mairesse~\cite{Mairesse}, the $A(n)$ being \mp  operators. The denomination rank~1 is natural for \mp operators.

We proved in~\cite{GM} that this property is generic for i.i.d. \mp operators: it is fulfilled when the support of the law of $A(1)$ is not included the union of finitely many affine hyperplanes.

Although this result could suggest the opposite, the MLP depends on the law of $A(1)$, and not only on its support :
if $\left(U(n)\right)_{n\in\N}$ is an i.i.d. sequence with the support of $U(1)$ equal to $[0,1]$, and $A(n)$ are the \mp operators defined by the matrices
$$A(n)=\left(
\begin{array}{cc}
-U(n) & 0 \\ 
0 & -U(n)
\end{array}\right),$$
then $\left(A(n)\right)_{n\in\N}$ has the MLP property iff $\P(U(n)=0)>0$.

The weaker condition that there is an operator with rank~1 in the closed semigroup generated by the support of  the law of $A(1)$ has been investigated by J. Mairesse for \mp operators. It ensures the weak convergence of $\overline{x}(n,.)$ but does not seem appropriate for our construction. 

\subsection{Known results}
Before describing our analysis, we give a brief review of published limit theorems about $x(n,X^0)$.

There has been many papers about the law of large numbers for products of random \mp matrices since it was introduced by J.E~Cohen~\cite{Cohen}. Let us cite F.~Baccelli~\cite{Baccelli}, the last one by T.~Bousch and J.~Mairesse~\cite{BouschMairesseEng} and our PhD thesis~\cite{theseGM} (in French). The last article proves results for a larger class of topical operators, called uniformly topical.

J.M. Vincent has proved a law of large number for topical operators, that will be enough in our case~:
\begin{thm}[\cite{vincent}]\label{thVincent}
Let $\left(A(n)\right)_{n\in\N}$ be a stationary ergodic sequence of topical operators and $X^0$ an  $\R^d$-valued random variable. If $A(1).0$ and $X^0$ are integrable, then there exists $\overline{\gamma}$ and $\underline{\gamma}$ in $\R$ such that 
\begin{eqnarray*}
\lim_n \frac{\max_ix_i(n,X^0)}{n}&=&\overline{\gamma}~\textrm{a.s.}\\
\lim_n \frac{\min_ix_i(n,X^0)}{n}&=&\underline{\gamma}~\textrm{a.s.}
\end{eqnarray*}
\end{thm}

F. Baccelli and J. Mairesse give a condition to ensure $\overline{\gamma}=\underline{\gamma}$, hence the convergence of $\frac{x(n,X^0)}{n}$:
\begin{thm}[\cite{BM96}]\label{LGN}
Let $\left(A(n)\right)_{n\in\N}$ be a stationary ergodic sequence of topical operators and $X^0$ an  $\R^d$-valued random variable such that $A(1).0$ and $X^0$ are integrable. If there exists an $N$, such that $A(N)\cdots A(1)$ has a bounded projective image with positive probability, then there exists $\gamma$ in $\R$ such that 
$$ \lim_n \frac{x(n,X^0)}{n}=\gamma\1~\textrm{a.s.}$$
\end{thm}
In this case $\gamma$ is called the Lyapunov exponent of the sequence. We notice that the MLP property implies a bounded projective image with positive probability.

The following result has been proved by J. Mairesse when the $A(n)$ are \mp operators, but can be extended to topical operators with the same proof. It will be the key point to ensure the spectral gap.
\begin{thm}[Mairesse \cite{Mairesse}]\label{strcoupling}
If the stationary and ergodic sequence $\left(A(n)\right)_{n\in\Z}$ of random variables with values in $Top_d$  has memory loss property, then there exists a random variable $Y$ with values in $\prm$ such that $Y_n:=A(n)\cdots A(1)Y$ is stationary. Moreover 
$$\lim_{n\rightarrow\infty}\P\left(\exists x_0, Y_n\neq \bar{x}(n,x_0)\right)=0.$$
In particular $\overline{x}(n,x_0)$ converges in total variation uniformly in $x_0$.
\end{thm}
The law of $Y$ is called the invariant probability measure.

To end this section we mention  two limit theorems, which are close to ours, but obtained by different ways. We will compare those results to ours in section~\ref{commentaires}.

With a martingale method J. Resing and al.~\cite{RVH} have obtained a central limit theorem for $x(n,X^0)$, when the Markov chain $\overline{x}(n,.)$ is aperiodic and uniformly $\Phi$-recurrent. The theorem has been stated for \mp operators, but it should make no difference to use topical ones. 

With a subadditivity method, F. Toomey~\cite{toomey} has proved a large deviation principle for $x(n,x_0)$ when the projective image of $A(N)\cdots A(1)$ is bounded.

\subsection{Principle of the analysis}\label{principes}
From now on, $\left(A(n)\right)_{n\in\N}$ is an i.i.d sequence of topical operators with the MLP property.

The first step of the proof is to split our Markov chain $x(n,.)$ into another Markov chain and a sum of cocycles over this chain, following what \'E. Le Page made for products of random matrices. For any  topical function $\phi$ from $\R^d$ to $\R$, $\phi(Ax)-\phi(x)$ only depends on $A$ and $\overline{x}$. Therefore $\phi(x(n,.))-\phi(x(n-1,.))$ only depends on $A(n)$ and $\overline{x}(n-1,.)$. Since $\prm$ can be seen as an hyperplane of $\R^d$, $x(n,.)$ can be replaced by $\left(\phi(x(n,.)),\overline{x}(n,.)\right)$. (cf. lemma~\ref{bilip})

According to theorem~\ref{strcoupling}, we know that $\overline{x}(n,.)$ converges. On the other hand, by theorem~\ref{LGN} $x(n,X^0)$, goes to infinity (if $\gamma\neq 0$) in the direction of $\1$, so $\phi(x(n,.))\sim \gamma n$. We investigate the oscillations of $\phi(x(n,.))-\gamma n $. Interesting $\phi$'s are defined by $\phi(x)=x_i$, $\phi(x)=\max_ix_i$, $\phi(x)=\min_ix_i$.\\

The second step is to prove the spectral gap for the operator defining the Markov chain $\left(A(n),\overline{x}(n-1,.)\right)_{n\in\N}$ and apply the results of~\cite{HH1} et~\cite{HH2} that give limit theorems for $\phi(x(n,X^0))-\phi(X^0)-\gamma n $. The spectral gap follows from the convergence of $\overline{x}(n,.)$, just like by \'E. Le Page~\cite{LePage}.

We use two series of results.
The first series are taken from the book by H. Hennion and L. Herv\'e~\cite{HH1} that sums up the classical spectral gap method developed since Nagaev~\cite{Nagaev} in a general framework. To apply it we demand integrability conditions on $\sup_x|\phi(A(1)x)-\phi(x)|$ to have a Doeblin operator on the space of bounded functions.
The second series are taken from the article~\cite{HH2} that is a new refinement of the method in the more precise framework of iterated Lipschitz operators. Since our model enters this framework, we get the same results with integrability conditions on $A(1)\,0$ that ensures that the Markov operator satisfy a Doeblin-Fortet condition on functions spaces defined by weights. The comparison between the two series of results will be made in section~\ref{commentaires}.

\section{Statement of the limit theorems}\label{statements}
\subsection{General case}
From now on, we state the results that we will prove in  section~\ref{proofs}.

For local limit theorem and for renewal theorem we need non arithmeticity conditions.
There are three kind of non arithmeticity, depending if the theorem follows from~\cite{HH1} or~\cite{HH2}. We will denote them respectively by (weak-) non arithmeticity and algebraic non arithmeticity.
When $d=1$ they  fall down to the usual non arithmeticity condition for real i.i.d. variables.
Algebraic non arithmeticity will be defined before the statement of LLT, but other non arithmeticity conditions will be defined in section~\ref{proofs} once we have given the definitions of the operator associated to the Markov chain. Unlike algebraic non arithmeticity, they depend on the 2-uple $\left(\left(A(n)\right)_{n\in\N},\phi\right)$, which will be called "the system".

Let $\left(A(n)\right)_{n\in\N}$  be an i.i.d sequence of topical operators with the MLP property.  The sequence $\left(x(n,.)\right)_{n\in\N}$ is defined by equation (\ref{defx}) and $\gamma$ is the Lyapunov exponent defined by theorem~\ref{LGN}.

Since the topology of the uniform convergence over compact subset on $Top_d$ has an enumerable basis of open sets, the  support of measures on it is well defined. We denote by $S_A$ the support of the law of $A(1)$ and by $T_A$ the  semi-group  generated by $S_A$ in $Top_d$.

\begin{thm}[CLT]\label{TCL}
Let $\left(A(n)\right)_{n\in\N}$  be an i.i.d sequence of topical operators with the MLP property and $X^0$ an $\R^d$-valued random variable independent from $\left(A(n)\right)_{n\in\N}$. Let $\phi$ be topical from $\R^d$ to $\R$.
Assume one of the following conditions:
\begin{enumerate}[i)]
\item $\sup_x\left|\phi(A(1)x)-\phi(x)\right|$ has a second moment,
\item $A(1)\,0$ has a $4+\epsilon$-th moment and $X^0$ has a  $2+\epsilon$-th moment.
\end{enumerate}
Then there exists $\sigma^2\ge0$ such that $\frac{x(n,X^0)-n\gamma\1}{\sqrt{n}}$ converges weakly to a random vector whose coordinates are equal and have law $\mathcal{N}(0,\sigma^2)$.\\

In the first case, or if $A(1)\,0$ has a $6+\epsilon$-th moment and $X^0$ has a $3+\epsilon$-th moment, then
\begin{itemize}
\item $\sigma^2=\lim\frac{1}{n}\E\left(\phi\left(x(n,X^0)\right)-n\gamma\right)^2 $
\item $\sigma=0$ iff there is a $\theta\in Top_d$ with rank~1  such that for any $A\in S_A$ and any $\theta'\in T_A$ with rank~1, $\theta A\theta'=\theta \theta' +\gamma\1$.
\end{itemize}
\end{thm}
\begin{rem}
According to lemma~\ref{invtheta}, if there is such a $\theta$, then every $\theta\in T_A$ with rank~1 has this property.
\end{rem}

\begin{thm}[CLT with rate]\label{TCLV}
Let $\left(A(n)\right)_{n\in\N}$  be an i.i.d sequence of topical operators with the MLP property and $X^0$ an $\R^d$-valued random variable independent from $\left(A(n)\right)_{n\in\N}$. Let $\phi$ be topical from $\R^d$ to $\R$.
Assume one of the following conditions:
\begin{enumerate}[i)]
\item $\sup_x\left|\phi(A(1)x)-\phi(x)\right|$ has an $l$-th moment with $l\ge 3$,
\item $A(1)\,0$ has an $l$-th moment, with $l>6$.
\end{enumerate} 
If $\sigma^2>0$ in theorem~\ref{TCL}, then there exists $C\ge 0$ such that for every initial condition $X^0$ with an $l$-th moment, we have
\begin{eqnarray}\label{vitTCL}
\lefteqn{\hspace{-3.5cm}\sup_{u\in\R}\left|\P[\phi(x(n,X^0))-n\gamma-\phi(X_i^0)\le\sigma u\sqrt{n}] -\mathcal{N}(0,1)(]-\infty,u])\right|}\nonumber\\
&&\le\frac{C\left(1+\E\left[\left\|X^0\right\|_\infty^l\right]\right)}{\sqrt{n}},
\end{eqnarray} 
\begin{eqnarray*}
\lefteqn{\hspace{-2cm}\sup_{u\in\R^d}\left|\P[x(n,X^0)-n\gamma\1\le\sigma u\sqrt{n}]-\mathcal{N}(0,1)(]-\infty,\min_iu_i])\right|}\\
&&\le \frac{C\left(1+\E\left[\left\|X^0\right\|_\infty^l\right]+\E\left[\left\|A(1)0\right\|_\infty^l\right]\right)}{n^{{\frac{l}{2(l+1)}}}}.\end{eqnarray*}
\end{thm}

\begin{defn}
We say that the sequence $\left(A(n)\right)_{n\in\N}$ is algebraically arithmetic if there are $a,b\in\R$ and a $\theta\in Top_d$ with rank~1  such that for any $A\in S_A$ and any $\theta'\in T_A$ with rank~1,
\begin{equation}\label{eqANA}
(\theta A\theta'-\theta \theta')(\R^d) \subset (a+b\Z)\1.
\end{equation}
Otherwise the sequence is algebraically non arithmetic.
\end{defn}
\begin{rem}

According to lemma~\ref{invtheta}, if there is such a $\theta$, then every $\theta\in T_A$ with rank~1 has this property.

Moreover, for any $\theta,\theta'\in Top_d$ with rank~1 and any $A\in Top_d$, the function $\theta A\theta'-\theta \theta'$ is constant with value in $\R\1$.
\end{rem}

\begin{thm}[LLT]\label{TLL}
Let $\left(A(n)\right)_{n\in\N}$  be an i.i.d sequence of topical operators with the MLP property and $X^0$ an $\R^d$-valued random variable independent from $\left(A(n)\right)_{n\in\N}$. Let $\phi$ be topical from $\R^d$ to $\R$.
Assume one of the following conditions:
\begin{enumerate}[i)]
\item $\sup_x\left|\phi(A(1)x)-\phi(x)\right|$ has a second moment, $\sigma>0$,  and the system is non arithmetic
\item $A(1)\,0$ has a $4+\epsilon$-th moment, $X^0\in \mathbb{L}^\infty$ and the sequence $\left(A(n)\right)_{n\in\N}$ is algebraically non arithmetic.
\end{enumerate}
Then $\sigma>0$ and there exists a $\sigma$-finite measure $\alpha$ on $\R^d$, so that for any continuous function $h$  with compact support, we have:
$$\lim_n\sup_{u\in\R}\left|\sigma\sqrt{2\pi n}\E\left[h\left(x(n,X^0)-n\gamma\1-u\1\right)\right]-\E\left[e^{-\frac{(u+\phi(X^0))^2}{2n\sigma^2}}\right]\alpha(h)\right|=0.$$
Moreover the image of $\alpha$  by  the function $x\mapsto (\overline{x},\phi(x))$ is the product of the invariant probability measure on $\prm$ by the Lebesgue measure.
\end{thm}

\begin{rem}
Like in the usual LLT, this theorem says that the probability for $x(n,X^0)$ to fall in a box  decreases like $\frac{1}{\sqrt{n}}$.
To replace the continuous functions by indicator functions of the box, we need to know more about the invariant probability measure on $\prm$. In particular, numerical simulations show that some hyperplanes may have a weight for this probability measure, so those hyperplanes could not intersect the boundary of the box.
\end{rem}

The algebraic non arithmeticity is optimal in the following sense:
\begin{prop}\label{optNA}
If the conclusion of theorem~\ref{TLL} is true, then $\left(A(n)\right)_{n\in\N}$ is algebraically non arithmetic.
\end{prop}

\begin{thm}[Renewal theorem]\label{renouv}
Let $\left(A(n)\right)_{n\in\N}$  be an i.i.d sequence of topical operators with the MLP property and $X^0$ an $\R^d$-valued random variable independent from $\left(A(n)\right)_{n\in\N}$.
Assume that there is a topical $\phi$ from $\R^d$ to $\R$ such that $\sup_x\left|\phi(A(1)x)-\phi(x)\right|$ has a second moment. We denote by $\alpha$ the same measure as in theorem~\ref{TLL}. If  $\gamma>0$ and  the system is weakly non arithmetic, then for any function $h$ continuous with compact support and any initial condition $X^0$, we have:
$$\lim_{a\rightarrow-\infty} \sum_{n\ge1}\E\left[h\left(x(n,X^0)-a\1\right)\right]=0,$$ 
$$\lim_{a\rightarrow+\infty} \sum_{n\ge1}\E\left[h\left(x(n,X^0)-a\1\right)\right]=\frac{\alpha(h)}{\gamma}.$$
\end{thm}

\begin{rem}
The vector $\1$ gives the average direction in which $x(n,X^0)$ is going to infinity. Like in the usual renewal theorem, this theorem says that the average number of $x(n,X^0)$ falling in a box is asymptotically proportional to the length of this box, when the box is going to infinity in that direction.
Like in the LLT, to replace the continuous functions by indicator functions of the box, we need to know more about the invariant probability measure on $\prm$.
\end{rem}

\begin{thm}[Large deviations]\label{PGD}
Let $\left(A(n)\right)_{n\in\N}$  be an i.i.d sequence of topical operators with the MLP property and $X^0$ an $\R^d$-valued random variable independent from $\left(A(n)\right)_{n\in\N}$. Let $\phi$ be topical from $\R^d$ to $\R$.
Assume that $\sup_x\left|\phi(A(1)x)-\phi(x)\right|$  has an exponential moment, and that $\sigma^2>0$ in theorem~\ref{TCL}. Then, there exists a non negative strictly convex function $c$, defined on a neighborhood of $0$ and vanishing only at $0$ such that for any bounded initial condition $X^0$ and any $\epsilon>0$ small enough we have:
$$\lim_{n}\frac{1}{n}\ln\left(\P\left[\phi\left(x(n,X^0)\right)-n\gamma>n\epsilon\right]\right)=-c(\epsilon),$$
$$\lim_{n}\frac{1}{n}\ln\left(\P\left[\phi\left(x(n,X^0)\right)-n\gamma<-n\epsilon\right]\right)=-c(-\epsilon).$$
\end{thm}

\subsection{Max-plus case}
When the $A(n)$ are \mp operators, it is natural to chose $\phi(x)=\max_ix_i$. In this case we get $\min_j\max_iA_{ij}\le\phi(Ax)-\phi(x)\le\max_{ij}A_{ij}$, so integrability condition can be checked on the last two quantities.

\begin{thm}[CLT]\label{TCL1mp}
Let $\left(A(n)\right)_{n\in\N}$  be an i.i.d sequence of \mp operators with the MLP property and $X^0$ an $\R^d$-valued random variable independent from $\left(A(n)\right)_{n\in\N}$.
If $\max_{ij}A(1)_{ij}$ and $\min_{j}\max_{i}A(1)_{ij}$ have a second moment, then there exists $\sigma^2\ge0$ such that for every initial condition $X^0$,
\begin{enumerate}[(i)]
\item $\frac{x(n,X^0)-n\gamma\1}{\sqrt{n}}$ converges weakly to a random variable whose coordinates are equals and have law $\mathcal{N}(0,\sigma^2)$,
\item $\sigma^2=\lim\frac{1}{n}\int \left(\max_{i,j}A(n)\cdots A(1)_{ij}\right)^2d\P$.
\end{enumerate}
\end{thm}
Theorems~\ref{TCLV}~to~\ref{PGD} are specialized in the same way: the conclusion is valid if $\phi(x)=\max_ix_i$ and the moment hypothesis on $\sup_x|\phi(A(1)x)-\phi(x)|$ is satisfied by $\max_{ij}A(1)_{ij}$ and $\min_{j}\max_{i}A(1)_{ij}$.

In this case, we also get another condition to avoid degeneracy in the CLT.
To state it, we recall a few definitions and results about \mp matrices:
\begin{defn}\label{defprod}\
For any $k,l,m\in\N$, the product of two matrices $A\in\rma^{k\times l}$ and $B\in\rma^{l\times m}$ is the matrix $A B\in\rma^{k\times m}$ defined by~:
$$\forall 1\le i\le k, \forall 1\le j\le m, (AB)_{ij}:=\max_{1\le p\le l}A_{ip}+ B_{pj}.$$
\end{defn}
If those matrices have no line of $-\infty$, then the \mp operator defined by $AB$ is the composition of those defined by $A$ and the one defined by $B$.

\begin{defn}
A circuit on a directed graph is a closed path on the graph.
Let $A$ be a square matrix of size $d$ with entries in $\rma$.
\begin{enumerate}[i)]
\item The graph  of $A$ is the directed weighted graph whose nodes are the integers from $1$ to $d$ and whose arcs are the $(i,j)$ such that $A_{ij}>-\infty$. The weight on $(i,j)$ is $A_{ij}$. The graph will be denoted by $\mathcal{G}(A)$ and the set of its elementary circuits  by $\mathcal{C}(A)$.
\item The average weight of a circuit $c=(i_1,\cdots,i_n,i_{n+1})$ (where $i_1=i_{n+1}$) is $aw(A,c):=\frac{1}{n}\sum_{j=1}^n A_{i_ji_{j+1}}.$
\item The \mp-spectral radius\footnote{this quantity is the maximal \mp-eigenvalue of $A$, that is $$\rm(A)=\max\{\lambda\in\rma|\exists V\in\rma^d\backslash\{(-\infty)^d\}, AV=V+\lambda\1\}.$$ See~\cite{theseGaubert}.} of $A$ is $\rm(A):=\max_{c\in\mathcal{C}(A)}aw(A,c)$.
\end{enumerate}
\end{defn}

\begin{thm}\label{s>0}
Assume the hypothesis of theorem~\ref{TCL1mp}, with $\gamma=0$. Then the variance $\sigma^2$ in theorem~\ref{TCL1mp} is $0$ if and only if $\left\{\rm(B)|B\in T_A\right\}=\{0\}$.
\end{thm}

Theorem 3.2 of~\cite{GM} gives a condition to ensure the memory loss property. This condition also ensures that there are two matrices in $S_A$ with  two distinct spectral radius. This proves the following corollary:
\begin{cor}\label{thgenesupp}
Let the law of $A(1)$ be a probability measure on the set of $d\times d$ matrices with finite second moment whose support is not included in the union of finitely many affine hyperplanes of $\R^{d\times d}$. Then $x(n,.)$ satisfies the conclusions of theorem~\ref{TCL1mp} with $\sigma >0$.
\end{cor}

We also give a sufficient condition to ensure the algebraic non arithmeticity:
\begin{thm}\label{NA}
Assume the hypothesis of theorem~\ref{TLL} $ii)$ except the algebraic non arithmeticity and $A(n)$ are \mp operators. If $\left(A(n)\right)_{n\in\N}$ is algebraically arithmetic, then there are $a,b\in\R$ such that $$\{\rm(B)|B\in S_A, \mathcal{G}(B) \textrm{strongly connected}\}\subset a+b\Z.$$
\end{thm}

Together with corollary~\ref{thgenesupp}, this proves that the hypothesis are generic in the following sense:
\begin{cor}\label{thgeneNA}
If the law of $A(1)$ is a probability measure on the set of $d\times d$ matrices with $4+\epsilon$-th moment whose support is not included in the union of enumerably many affine hyperplanes of $\R^{d\times d}$, then $x(n,.)$ satisfies the conclusions of theorem~\ref{TLL}.
\end{cor}

\subsection{Comments}\label{commentaires}
The following table sums up the limit theorems. In each situation we assume that the sequence $\left(A(n)\right)_{n\in\N}$ has the memory loss property.
\begin{center}\begin{tabular}{|c|c|c|c|}
\hline Theorems: &\multicolumn{2}{|c|}{Moments of}&Additional \\
&$A(1)\,0$&$\max_{ij}A(1)_{ij}$ and $\min_{j}\max_{i}A(1)_{ij}$& condition\\ 
\hline CLT & $4+\epsilon$ & $2$ &\\ 
\hline CLT with rate & $6+\epsilon$ & $3$ & \\ 
\hline LLT & $4+\epsilon$ & $2$ &NA\\
\hline Renewal& -- &$2$ &NA\\
\hline LDP & -- & exp & $X^0\in L^\infty$\\ 
\hline \multicolumn{4}{c}{NA= non arithmeticity}\\
\end{tabular}\\
\end{center}

Let us first notice that the results of the second column are optimal in the sense that their restriction to $d=1$ are exactly the usual theorems for sum of i.i.d. real variable (except for LDP).

The results of the second column are stated for \mp operators because the $\sup_x\left|\phi(A(1)x)-\phi(x)\right|$ is not bounded for general topical operators. For other subclasses of topical operators, one has to choose $\phi$ such that $\sup_x\left|\phi(A(1)x)-\phi(x)\right|$ is integrable. For instance $\phi(x)=\min_ix_i$ is natural for $(\min,+)$ operators. Actually, it should be possible to derive renewal theorem and large deviation principle with the method of~\cite{HH2} but this has not been written down.

The results of the first column require stronger integrability conditions but they are also better for two reasons: they are true for any topical operators and the algebraic non arithmeticity does not depend on $\phi$. It is expressed without introducing the Markovian operator $Q$ although the system is algebraically non arithmetic iff $Q$ has an eigenvector with eigenvalue with modulus~$1$. Moreover for \mp operators the algebraic non arithmeticity can be deduced from theorem~\ref{NA}. \\

An important case for \mp operators is when $A_{ij}\in\R^+\cup\{-\infty\}$ and $A_{ii}\ge 0$ because it modelizes situations where $x_i(n,.)$ is the date of the $n$-th event of type $i$ and the $A_{ij}$ are delays. In this case the integrability of $\max_{ij}A(1)_{ij}$ and $\min_{j}\max_{i}A(1)_{ij}$ is equivalent to the integrability of $A(1)\,0$.\\

We mentioned earlier that J. Resing and al.~\cite{RVH} obtained a central limit theorem. In a sense our result is weaker because the MLP property implies that  $\overline{x}(n,.)$ is uniformly $\Phi$-recurrent and aperiodic. But our integrability conditions are much weaker and the MLP property is easier to check.

F. Toomey's large deviations principle only requires the uniform bound of the projective image that is a very strong  integrability condition. It suggests that the MLP property should not be necessary. But his formulation of the LDP is not equivalent to ours and in the \mp case it needs the fixed structure property, that is $\P(A_{ij}(1)=-\infty)\in\{0,1\}$.

\section{Proofs of the limit theorems}\label{proofs}
\subsection{From iterated topical functions to Markov chains}
In this section, we show that the hypothesis of the theorems on our model stated in section~\ref{statements} imply the hypothesis of the general theorems of~\cite{HH1} et~\cite{HH2}.
To apply the results of~\cite{HH1},
\begin{itemize}
\item the space $E$ will be $Top_d\times\prm$ with the Borel $\sigma$-algebra,
\item the transition probability $Q$ will be defined by $$Q\left((A,\bar{x}),D\right)=\P\left((A(1),\overline{Ax})\in D\right),$$
\item the Markov chain $X_n$ will be $\left(A(n),\bar{x}(n-1,.)\right)$,
\item the function $\xi$ will be defined by $\xi(A,\bar{x})=\phi(Ax)-\phi(x)$, where $\phi$ is a topical function from $\R^d$ to $\R$.
\item $\sigma(A)=\sup_x|\xi(A,x)|=\sup_x\left|\phi(Ax)-\phi(x)\right|<\infty$ a.s.\,.
\end{itemize}

With these definitions, $S_n=\sum_{l=1}^n\xi(X_l)$ is equal to $\phi(x(n,X^0))-\phi(X^0)$.

To apply the results of~\cite{HH1}, we still need to define the space $\B$.
\begin{defn}
Let $\mathcal{L}^\infty$ be the space of complex valued bounded continuous functions on $\prm$.

Let $j$ be the function from $Top_d\times\prm$ to $\prm$ such that $j(A,\bar{x})=\overline{A  x}$ and $I$ the function from $\R^{\prm}$ to $\R^{(Top_d\times\prm)}$  defined by
$$I(\phi)=\phi\circ j.$$
We call $\B^\infty$ the image of $\mathcal{L}^\infty$ by $I$.
\end{defn}

Since $(\mathcal{L}^\infty,\|.\|_\infty)$ is a Banach space,  $I$ is an injection, and $\|\phi\circ I\|_\infty =\|\phi\|_\infty$, $(\B^\infty,\|.\|_\infty)$ is also a Banach space.

\begin{defn}
The Fourier kernels denoted by $Q_t$ or $Q(t)$ are defined for any $t\in \C$ by
$$Q_t(x,dy)=e^{it\xi(y)}Q(x,dy).$$
We say that the system is non arithmetic if $Q(.)$ is continuous from $\R$ to the space $\mathcal{L}_{\B^\infty}$ of continuous linear operators on ${\B^\infty}$ and, for any $t\in\R^*$, the spectral radius of~$Q(t)$ is strictly less than~$1$.

We say that the system is weakly non arithmetic if $Q(.)$ is continuous from $\R$ to $\mathcal{L}_{\B^\infty}$ and, for all $t\in\R^*$, $Id-Q(t)$, where $Id$ is the identity on ${\B^\infty}$, is invertible.
\end{defn}

\begin{prop}\label{Hm}
If $\sigma\left(A(1)\right)$ has an $m$-th moment and if $A(n)$ has the MLP property, then $(Q,\xi,\B^\infty)$ satisfies condition $\mathcal{H}(m)$ of~\cite{HH1}.
Moreover, the interval $I_0$ in condition $(H3)$ is the whole $\R$ and $s(Q,\B^\infty)=1$.
\end{prop}

To prove $(H2)$, we will use theorem~\ref{strcoupling}.
In the sequel, $\nu_0$ will be the law of $Y$ in that theorem, that is the (unique) invariant probability measure.

\begin{proof}[Proof of proposition~\ref{Hm}]
Condition $(H1)$ is trivial, because of the choice of~$\B^\infty$.\\

To check condition $(H2)$, we take $\nu:=\mu\otimes \nu_0$. It is $Q$-invariant by definition of $\nu_0$. This proves $(i)$. To prove $(ii)$ and $(iii)$, we investigate the iterates of $Q$. For any $\phi\in\mathcal{L}^\infty$, and $x\in\R^d$ we have:
\begin{eqnarray*}
\left| Q^n(\phi\circ j)(A,\bar{x})-\nu(\phi\circ j)\right|
&=& \left|Q^n(\phi)(\overline{A  x})-\nu_0(\phi)\right|\\
&=&\left|\int \phi\left(\bar{x}(n,\overline{A  x})\right) d\P-\int \phi(Y_n)d\P\right|\\
&\le &\|\phi\|_\infty 2\P\left(\exists x_0, Y_n\neq \bar{x}(n,x_0)\right),
\end{eqnarray*}

If we denote $\psi\mapsto \nu(\psi)$ by $N$, we obtain:
\begin{equation}
\|Q^n-N\|\le 2\P\left(\exists x, Y_n\neq \bar{x}(n,x)\right)\rightarrow 0.
\end{equation} 
This proves that the spectral radius $r(Q_{|Ker N})$ is strictly less than~$1$ and that $\sup_n\|Q^n\|<\infty$. Since $Q_{|Im N}$ is the identity, $dim (Im N)=1$ and $\B^\infty=Ker N\oplus Im N$, $Q$ is quasi-compact, so $(ii)$ is checked. Moreover $s(Q,\B^\infty)=1$.
It also proves $Ker(1-Q)\subset Im N$, which implies $(iii)$.\\

To prove $(H3)$ we set $Q_t^{(k)}:=e^{it\xi(y)}(i\xi(y))^kQ(x,dy)$ and 
\begin{equation}
\Delta_h^{(k)}:= Q^{(k)}_{t+h}-Q^{(k)}_t-hQ^{(k+1)}_t.
\end{equation}
To prove that $Q^{k+1}$ is the derivative of $Q^{k}$, it remains to bound $\|\frac{1}{h}\Delta_h^{(k)}\|$ by a quantity that tends to zero with $h$.

To this aim, we introduce the following function:
$$\left\{\begin{array}{lccl}f:&\R &\rightarrow &\C\\&t&\mapsto &e^{it}-1-it.\end{array}\right.$$
The calculus will be based on the  following estimations on $f$: $|f(t)|\le 2t$, and $|f(t)|\le t^2$.

Now everything follows from  
\begin{equation}\label{introf}
\Delta_h^{(k)}(\phi\circ j)(A,\bar{x})
= \int  \phi\left(\overline{BAx}\right)e^{it\xi(B,\overline{Ax})}\left(i\xi(B,\overline{A x})\right)^k f\left(h\xi(B,\overline{Ax})\right) d\mu(B) .
\end{equation}

First it implies that
\begin{equation}\label{norminfinie}
\left\|\Delta_h^{(k)}(\phi\circ j)\right\|_\infty
\le \|\phi\|_\infty  \int \sigma^k(B)\left\|f\left(h\xi(B,.)\right) \right\|_\infty d\mu(B).
\end{equation}

Since $|f(t)|\le t^2$,
$$\frac{1}{|h|}\sigma(B)^{k} \left\| f\left(h\xi(B,.)\right) \right\|_\infty    \le h   \sigma^{k+2}(B) \rightarrow 0 .$$
Since $|f(t)|\le 2t$, 
$$\frac{1}{|h|}\sigma^{k}(B) \left\| f\left(h\xi(B,.)\right) \right\|_\infty    \le 2  \sigma^{k+1}(B).$$

When $k<m$, $\sigma^{k+1}$ is integrable so the dominate convergence theorem and the last two equations show that
\begin{equation}\label{term1}
\int\sigma(B)^{k} \left\| f\left(h\xi(B,.)\right) \right\|_\infty d\mu(B)=o(h).
\end{equation} 

Finally for any $k<m$, $\|\frac{1}{h}\Delta_h^{(k)}\|$ tends to zero, so $Q_t^{(k+1)}$ is the derivative of $Q_.^{(k)}$ in $t$.

To prove that $Q_.^{(m)}$ is continuous, we notice that
\begin{equation}
\left(Q_{t+h}^{(m)}-Q_t^{(m)}\right)(\phi\circ j)(A,\overline{x})= \int  \phi\left(\overline{BAx}\right)e^{it\xi(B,\overline{Ax})}\left(i\xi(B,\overline{A x})\right)^m g\left(h\xi(B,\overline{Ax})\right) d\mu(B) .
\end{equation}
where $g(t)=e^{it}-1$. Then we apply the same method as before, replacing the estimates on $f$ by $|g(t)|\le t$ to prove the convergence, and by  $|g(t)|\le 2$ to prove the domination.

This proves $(H3)$ and the additional assumption of proposition~\ref{Hm}.
\end{proof}

In their article~\cite{HH2} H. Hennion and L. Herv\'e have proved limit theorems for sequences $\xi(Y_n,Z_{n-1})$, where $(Y_n)_{n\in\N}$ is an i.i.d. sequence of Lipschitz operators on a metric space $\M$, and $Z_n$ is defined by $Z_{n+1}=Y_{n+1}Z_n$. As explained in section~\ref{principes}, we take $\M=\prm$, $Y_n=A(n)$ and again $\xi(A,\overline{x})=\phi(Ax)-\phi(x)$. In this case $Z_n=\overline{x}(n,X^0)$ and $S_n=\phi\left(x(n,X^0)\right)-\phi(X^0)$. Moreover in our situation, the $Y_n$, which are the projective function defined by $A(n)$, are 1-Lipschitz. Following the same proof as~\cite{HH2} with this additional condition, we get the CLT (resp. CLT with rate, LLT) for $S_n$ under the hypothesis of theorem~\ref{TCL} (resp.~\ref{TCLV},\ref{TLL}) on $A(1)0$.

The integrability conditions are weaker than in~\cite{HH2}, because the Lipschitz coefficient is uniformly bounded. The only difference in the proof is the H\"older inequality of the 4th part of proposition~7.3 of~\cite{HH2}: the exponents in the inequality should be changed to $1$ and $\infty$.\\
 
Let us give the notations of~\cite{HH2} we need to state the results. $G$ is the semi-group of the operators on $\M$, 1-Lipschitz for distance $\delta$. For a fixed $x_0\in\M$, every $\eta\ge 1$ and every $n \in \N$, we set $\M_\eta=\E[\delta^\eta(Y_1x_0,x_0)]$ and $\Cr_n=\E[c(Y_n\cdots Y_1)]$, where $c(.)$ is the Lipschitz coefficient.

When there is an $N\in\N$ such that $\Cr_N<1$, there is a $\lambda_0 \in]0,1[$, such that \mbox{$\int_Gc(g)\left(1+\lambda_0 \delta(gx_0,x_0)\right)^{2\eta}d\mu^{*N}(g)<1$.} We chose one such $\lambda_0$ and set the following notations:
\begin{enumerate}[(i)]
\item $\B_\eta$ is the set of functions $f$ from $\M$ to $\C$ such that $m_\eta(f)<\infty$, with norm $\| f\|_\eta=|f|_\eta+m_\eta(f)$, where
\begin{eqnarray*}
|f|_\eta&=&\sup_x\frac{|f(x)|}{(1+\lambda_0 \delta(x,x_0))^{1+\eta}}, \\
m_\eta(f)&=&\sup_{x\neq y}\frac{|f(x)-f(y)|}{\delta(x,y)\left(1+\lambda_0\delta(x,x_0)\right)^\eta\left(1+\lambda_0 \delta(y,x_0)\right)^\eta}.
\end{eqnarray*}
\item We say that the system is $\eta$-non arithmetic if there is no $t\in\R\backslash\{0\}$, no $\rho\in\C$, and no $w\in\B_\eta$ with non-zero constant modulus on the support $S_{\nu_0}$ of the invariant probability measure $\nu_0$ such that $|\rho|=1$ and for all $n\in\N$, we have
\begin{equation}\label{eqNA}
e^{itS_n}w(Z_n)=\rho^nw(Z_0) \P-\textrm{ a.s.,}
\end{equation}
when $Z_0$ has law $\nu_0$.
\end{enumerate}

\begin{rem}[non arithmeticity]
In the first frame the non arithmeticity condition is about the spectral radius of $Q_t$. Here we work with the associated $P_t$ that acts on $\M$ instead of $G\times\M$ (cf.~\cite{HH2}). If $P_t$ is quasi-compact, then the spectral radius $r(P_t)$ is~$1$ iff $P_t$ has an eigenvalue $\rho$ with modulus~$1$.
It is shown in proposition~9.1'~of~\cite{HH2} that if $r(P_t)=1$, then $P_t$ is quasi-compact as an operator on $\B_\eta$ and that an eigenvector $w$ with eigenvalue $\rho$ satisfies equation~(\ref{eqNA}).
\end{rem}

\begin{prop}\label{HS}\ 
\begin{enumerate}
\item If $A(1)\,0$ has an $\eta$-th moment, with $\eta\in\R^+$, then $\M_\eta<\infty$. If the sequence has the MLP property, then there is $n_0\in\N$ such that $\mathcal{C}_{n_0}<1$. If $\overline{X^0}$ has an $\eta$-th moment $\eta\in\R^+$, then $f\mapsto\E[f(\overline{X^0})]$ is continuous on $B_\eta$.
\item Algebraic non arithmeticity implies $\eta$-non arithmeticity for any $\eta>0$.
\end{enumerate} 
\end{prop}

The first part of the proposition is obvious. The second part relies on the  next two lemma that will be proved after the proposition:
\begin{lem}\label{suppnu}
The support of the invariant measure $\nu_0$ is $$S_{\nu_0}:=\overline{\{\overline{\theta\1}|\theta\in T_A, \theta\textrm{ with rank~1}\}}.$$
\end{lem}

\begin{lem}\label{invtheta}
If equation~(\ref{eqANA}) is satisfied by some $\theta$ with rank~1, any $A\in S_A$ and any $\theta'\in T_A$ with rank~1, it is satisfied by any $\theta \in T_A$ with rank~1.
\end{lem}

\begin{proof}[Proof of proposition~\ref{HS}]
Let us assume that the system is $\eta$-arithmetic. Then there are $w\in\B_\eta$ and $t,a\in\R$ such that for $\mu$-almost every $A$ and $\nu_0$ almost every $\overline{x}$, we have:
\begin{equation}\label{eqNA1}
e^{it\left(\phi(Ax)-\phi(x)\right)}w(\overline{Ax})=e^{ita}w(\overline{x}).
\end{equation}
Since all functions in this equation are continuous, it is true for $\overline{x}\in S_{\nu_0}$ and $A\in T_A$. Since $S_{\nu_0}$ is $T_A$ invariant, we iterate equation~(\ref{eqNA1}) and get 
\begin{equation}\label{eqNAn}
e^{it\left(\phi(Tx)-\phi(x)\right)}w(\overline{Tx})=e^{itan_T}w(\overline{x}),
\end{equation}
where $T\in T_A$ and $n_T$ is the number of operators of $S_A$ to be composed to obtain $T$.

Because of the MLP property, there is a $\theta\in T_A$ with rank~1. For any $A\in S_A$, $\theta A\in T_A$, so we apply equation~(\ref{eqNAn}) for $T=\theta A$ and $T=\theta$ and divide the first equation by the second one. Since $n_{\theta A}=n_{\theta}+1$ and $\overline{\theta Ax}=\overline{\theta x}$ , we get
$$e^{it\left(\phi(\theta Ax)-\phi(\theta x)\right)}=e^{ita}.$$

Setting $b=\frac{2\pi}{t}$, it means that $\phi(\theta Ax)-\phi(\theta x)\in a+b\Z$. Since $\theta$ has rank one, $(\theta Ax-\theta x)\in\R\1$, so $\theta Ax-\theta x \in (a+b\Z)\1$, and the algebraic arithmeticity follows by lemma~\ref{suppnu}.
\end{proof}

\begin{proof}[Proof of lemma~\ref{suppnu}]

By theorem $\ref{strcoupling}$, there is sequence of random variables $Y_n$ with law $\nu_0$, such that $Y_n=A(n)\cdots A(1)Y$. Let $K$ be a compact subset of $\rma$ such that $Y\in K$ with positive probability.

For any $\theta\in T_A$ and any $\epsilon>0$, the set $V$ of topical functions $A$ such that $\delta(\overline{Ax},\overline{\theta x})\le\epsilon$ for all $\overline{x}\in K$ is a neighborhood of $\theta$. Therefore the probability for $A(n_\theta)\cdots A(1)$ to be in $V$ is positive and by independence of $Y$, we have:
$$\P\left[Y\in K ,A(n_\theta)\cdots A(1)\in V\right]>0.$$
Since $\overline{\theta\1}=\overline{\theta Y}$, this means that with positive probability, $$\delta\left(Y_{n_\theta},\overline{\theta \1}\right)=\delta\left(A(n_\theta)\cdots A(1)Y,\overline{\theta Y}\right)\le\epsilon,$$ so $\overline{\theta \1}\in S_ {\nu_0}.$

This proves that $\overline{\{\overline{\theta\1}|\theta\in T_A, \theta\textrm{ with rank~1}\}}\subset S_ {\nu_0} $.\\

In~\cite{Mairesse}, $\nu_0$ is obtained as the law of $Z=\lim_n\overline{A(1)\cdots A(n)\1}$. Indeed, the MLP property and the Poincar\'e recurrence theorem ensure that there are almost surely $M$ and $N$ such that $A(N)\cdots A(N+M)$ has rank~1. Therefore, for $n\ge N+M$, $\overline{A(1)\cdots A(n)\1}=\overline{A(1)\cdots A(N+M)\1}=Z$ . But $A(1)\cdots A(N+M)\in T_A$ almost surely, so \mbox{$Z\in\{\overline{\theta\1}|\theta\in T_A, \theta\textrm{ with rank~1}\}$} almost surely and $S_ {\nu_0}\subset \overline{\{\overline{\theta\1}|\theta\in T_A, \theta\textrm{ with rank~1}\}}$.
\end{proof}

\begin{proof}[Proof of lemma~\ref{invtheta}]
We assume that equation (\ref{eqANA}) is satisfied by $\theta=\theta_1$, any $A\in S_A$ and any $\theta'\in T_A$ with rank~1.

Let $A_1,\cdots, A_n \in S_A$, such that $\theta_2=A_1\cdots A_n$ has rank~1. For any $i\le n$, $A_i\cdots A_n\theta'$ has rank~1, so $(\theta_1A_{i}\cdots A_n\theta'-\theta_1A_{i+1}\cdots A_n\theta')(\R^d)\subset (a+b\Z)\1$.

Summing these inclusions for $i=1$ to $i=n$, we get $(\theta_1\theta_2\theta'-\theta_1\theta')(\R^d)\subset (na+b\Z)\1$ and 
\begin{equation}\label{eq1}
\left((\theta_1\theta_2A\theta'-\theta_1A\theta')-(\theta_1\theta_2\theta'-\theta_1\theta')\right)(\R^d)\subset b\Z\1.
\end{equation} 

Now we write $\theta_2\theta'$ as 
$$\theta_2\theta'=\theta_1\theta' +(\theta_1\theta_2\theta'-\theta_1\theta')- (\theta_1\theta_2\theta'-\theta_2\theta').$$

The last part does not depend on $\theta'$, so replacing $\theta'$ by $A\theta'$ and subtracting the first version, we get:
$$\theta_2A\theta'-\theta_2\theta'=\theta_1A\theta'-\theta_1\theta' +\left((\theta_1\theta_2A\theta'-\theta_1A\theta')-(\theta_1\theta_2\theta'-\theta_1\theta')\right).$$
With equation~(\ref{eq1}), this proves equation~(\ref{eqANA}) for $\theta=\theta_2$.

\end{proof}

\subsection{From Markov chains to iterated topical functions} 
Propositions~\ref{Hm} and~\ref{HS} prove that under the hypothesis of section~\ref{statements} the conclusions of the theorems of~\cite{HH1} and~\cite{HH2} are true. This gives results about the convergence of $\left(\phi\left(x(n,X^0)\right)-\phi\left(X^0\right)-n\nu(\xi),\overline{x}(n,X^0)\right)$.

When $\overline{X^0}$ has law $\nu_0$, the sequence $\left(A(n),\overline{x}(n,X^0)\right)_{n\in\N}$ is stationary, so it follows from Birkhoff theorem that $\gamma=\int\xi(A,\overline{x}) d\nu_0(\overline{x})d\mu(A)=\nu(\xi)$.

The following lemma will be useful to go back to $x(n,.)$.
\begin{lem}\label{bilip}
If $\phi$ is a topical function from $\R^d$ to $\R$, the function $\psi:x\mapsto (\phi(x),\overline{x})$ is a Lipschitz homeomorphism with Lipschitz inverse from $\R^d$ onto $\R\times\prm$.
\end{lem}
\begin{proof}
Let $(t,\overline{x})$ be an element of $\R\times\prm$. Then $\psi(y)=(t,\overline{x})$ if and only if there is an $a\in\R$ such that $y=x+a\1$ and $\phi(x)+a=t$. So the equation has exactly one solution $y=x+(t-\phi(x))\1$ and $\psi$ is invertible. 

It is well known that topical functions are Lipschitz, and the projection is linear, so it is Lipschitz and so is $\psi$.

For any $x,y\in\R^d$, we have $x\le y+\max_i(x_i-y_i)\1$, so $\phi(x)-\phi(y)\le \max_i(x_i-y_i)$. Therefore, for any $1\le i\le d$, we have 
$$\phi(x)-\phi(y)-(x_i-y_i)\le \max_i(x_i-y_i)-\min_i(x_i-y_i)=\delta(\overline{x},\overline{y}).$$ 
Permuting $x$ an $y$, we see that:
\begin{equation}\label{difftop}
|\phi(x)-\phi(y)-(x_i-y_i)|\le \delta(\overline{x},\overline{y}).
\end{equation}
Therefore $|x_i-y_i|\le |\phi(x)-\phi(y)|+\delta(\overline{x},\overline{y})$ and  $\psi^{-1}$ is Lipschitz.
\end{proof}

\begin{proof}[Proof of theorem~\ref{TCL}]
Without lost of generality, we assume that $\gamma=0$.
Theorem~A of~\cite{HH1} and proposition~\ref{Hm} or theorem~A of~\cite{HH2} and proposition~\ref{HS} prove that 
$\frac{\phi(x(n,X^0))-\phi(X^0)}{\sqrt{n}}$ converges to $\mathcal{N}(0,\sigma^2)$, which means that
$\frac{\phi(x(n,X^0))-\phi(X^0)}{\sqrt{n}}~\1$ converges to the limit specified in theorem~\ref{TCL}\,. We just estimate the difference between the converging sequence and the one we want to converge:
\begin{equation}\label{majdiff}
\Delta_n:=\left\|\frac{x(n,X^0)}{\sqrt{n}}-\frac{\phi\left(x(n,X^0)\right)-\phi(X^0)}{\sqrt{n}}~\1\right\|_\infty\le \frac{\left|\phi(X^0)\right|}{\sqrt{n}}+\frac{|\overline{x}(n,X^0)|_\mathcal{P}}{\sqrt{n}}.
\end{equation}

Each term of the last sum is a weakly converging sequence divided by $\sqrt{n}$ so it converges to zero in probability. This proves that $\Delta_n$ converges to zero in probability, which ensures the convergence of $\frac{x(n,X^0)}{\sqrt{n}}$ to the Gaussian law.\\

The expression of $\sigma^2$ is the direct consequence of theorems~A of~\cite{HH1} or theorem~S of~\cite{HH2}.\\

If $\sigma=0$, then again by theorem~A of~\cite{HH1}~or~S of~\cite{HH2}, there is a continuous function $\xi$ on $\prm$ such that 
\begin{equation}\label{eq2}
\phi(Ax)-\phi(x)=\xi(\overline{x})-\xi(\overline{Ax})
\end{equation} for $\mu$-almost every $A$ and $\nu_0$-almost every $\overline{x}$. 
Since all functions are continuous in this equation, (\ref{eq2}) is true for every $A\in S_A$ and $\overline{x}\in S_{\nu_0}$.
By induction we get it for $A\in T_A$ and if $\theta\in T_A$ has rank~1 and $\overline{x}\in S_{\nu_0}$, $\overline{\theta Ax}=\overline{\theta x}$, so $\phi(\theta Ax)=\phi(\theta x)$.

Since $\theta Ax-\theta x\in\R\1$, this means that $\theta Ax=\theta x.$
By lemma~\ref{suppnu}, it proves that $\theta A\theta'=\theta \theta'$ for any $\theta,\theta'\in T_A$ with rank~1 and $A\in S_A$.\\

Conversely, let us assume there is $\theta$ with rank one such that for any $\theta'\in T_A$ with rank~1 and $A\in S_A$, we have:
\begin{equation}\label{eq3}
\theta A\theta'=\theta \theta'.
\end{equation}
By lemma~\ref{invtheta} applied with $a=b=0$, it is true for any $\theta,\theta'\in T_A$ with rank~1, and any $A\in S_A$ and by induction, equation~(\ref{eq3}) is still true for $A\in T_A$.

Therefore, for any $m\in\N$ and $n\ge m+1$ and any $\theta'\in T_A$ with rank~1, if $A(n)\cdots A(n-m+1)$ has rank~1 ,then $x(n,\theta'\1)=A(n)\cdots A(n-m+1)\theta'\1$ and for any $N\in\N$
\begin{eqnarray}\label{eqxborne}
\lefteqn{\P\left( \|x(n,\theta'\1)\|_\infty\le N\right)}\nonumber\\
&\ge & \P\left( A(n)\cdots A(n-m+1)\textrm{has rank~1},  \|A(n)\cdots A(n-m+1)\theta'\1 \|_\infty\le N\right)\nonumber\\
&\ge & \P\left( A(m)\cdots A(1)\textrm{has rank~1},  \|A(m)\cdots A(1)\theta'\1 \|_\infty\le N\right).
\end{eqnarray}

We fix a $\theta'\in T_A$ with rank one. The MLP property says there is an $m$ such that $\P(A(m)\cdots A(1) \textrm{has rank~1})>0$. Therefore, there is an $N\in\N$ such that the right member of~(\ref{eqxborne}) is a positive number we denote by $\beta$.

Equation~(\ref{eqxborne}) now implies that for any $\epsilon>0$, if $n\ge \max(m,N^2\epsilon^{-2})$, then $\P( \|\frac{1}{\sqrt{n}}x(n,\theta'\1)\|_\infty\le \epsilon )\ge \beta$, so $\mathcal{N}(0,\sigma^2)[-\epsilon,\epsilon]\ge\beta$. When $\epsilon$ tends to zero, we get that $\mathcal{N}(0,\sigma^2)(\{0\})\ge\beta>0$, which is true only if $\sigma=0$.
\end{proof}

\begin{proof}[Proof of theorem~\ref{TCLV}]
Without loss of generality, we assume that $\gamma=0$.
Equation~(\ref{vitTCL}) follows from theorem~B of~\cite{HH1} and proposition~\ref{Hm} or from theorem~B of~\cite{HH2} and proposition~\ref{HS}

The only fact to check is that the initial condition defines a continuous linear form on $\B_\eta$, with norm at most $C\left(1+\E(\|X^0\|^l_\infty\right)$, that is for any $f\in\B_\eta$, we have:
$$|\E(f(X^0))|\le C\left(1+\E(\|X^0\|^l_\infty)\right)\|f\|_\eta.$$
It easily follows from the fact that $|f(x)|\le \|f\|_\eta (1+|x|_\mathcal{P})^{1+\eta}$ and $1+\eta \le l $.

Taking $y=0$ in~(\ref{difftop}), we get $|\phi(x)-x_i|\le|x|_\mathcal{P}$. Together with~(\ref{majdiff}) it proves that for any $u\in\R^d$ and any $\epsilon>0$
\begin{eqnarray}\label{majvit}
\lefteqn{ \P[x(n,X^0)\le\sigma u\sqrt{n}]}\nonumber\\
&\le & \P\left[\min_i x_i(n,X^0)\le\sigma \min_iu_i\sqrt{n}\right]\nonumber\\
&\le &\P\left[\phi(x(n,X^0))\le(\sigma\min_iu_i+2\epsilon)\sqrt{n}\right]
+ \P\left[ \frac{\left|\phi(X^0)\right|}{\sqrt{n}}\ge\epsilon\right]
+ \P\left[\frac{|\overline{x}(n,X^0)|_\mathcal{P}}{\sqrt{n}}\ge\epsilon\right]\nonumber\\
&\le&\mathcal{N}(0,1)(]-\infty,\min_iu_i+\frac{2\epsilon}{\sigma}])+\frac{C}{\sqrt{n}}
+\frac{\E\left(\left|\phi(X^0)\right|^l\right)}{(\epsilon\sqrt{n})^{l}}
+\frac{\E\left(|\overline{x}(n,X^0)|^l_\mathcal{P}\right)}{(\epsilon\sqrt{n})^{l}}\nonumber\\
&\le&\mathcal{N}(0,1)(]-\infty,\min_iu_i])+\frac{C}{\sqrt{n}}+\frac{2\epsilon}{\sigma}
+\frac{\E\left(\left|\phi(X^0)\right|^l\right)}{(\epsilon\sqrt{n})^{l}}
+\frac{\E\left(|\overline{x}(n,X^0)|^l_\mathcal{P}\right)}{(\epsilon\sqrt{n})^{l}}.
\end{eqnarray}

Conversely,
\begin{eqnarray}\label{minvit}
\lefteqn{ \P[x(n,X^0)\le\sigma u\sqrt{n}]}\nonumber\\
&\ge & \P\left[\phi(x(n,X^0))\le\sigma \min_iu_i\sqrt{n}\right]\nonumber\\
&\ge &\P\left[\phi(x(n,X^0))\le(\sigma\min_iu_i-2\epsilon)\sqrt{n}\right]
- \P\left[ \frac{\left|\phi(X^0)  \right|}{\sqrt{n}}\ge\epsilon\right]
- \P\left[\frac{|\overline{x}(n,X^0)|_\mathcal{P}}{\sqrt{n}}\ge\epsilon\right]\nonumber\\
&\ge &\mathcal{N}(0,1)(]-\infty,\min_iu_i-\frac{2\epsilon}{\sigma}])-\frac{C}{\sqrt{n}}
-\frac{\E\left(\left|\phi(X^0)  \right|^l\right)}{(\epsilon\sqrt{n})^{l}}
- \frac{\E\left(|\overline{x}(n,X^0)|^l_\mathcal{P}\right)}{(\epsilon\sqrt{n})^{l}}\nonumber\\
&\ge &\mathcal{N}(0,1)(]-\infty,\min_iu_i])-\frac{C}{\sqrt{n}}-\frac{2\epsilon}{\sigma}
-\frac{\E\left(\left| \phi(X^0) \right|^l\right)}{(\epsilon\sqrt{n})^{l}}
- \frac{\E\left(|\overline{x}(n,X^0)|^l_\mathcal{P}\right)}{(\epsilon\sqrt{n})^{l}}\nonumber\\.
\end{eqnarray}

Taking $\epsilon=n^{-{\frac{l}{2(l+1)}}}$ in (\ref{majvit}) and (\ref{minvit}) will conclude the proof of theorem~\ref{TCLV} if we can show that $\E\left(|\overline{x}(n,X^0)|^l_\mathcal{P}\right)$ is bounded uniformly in $n$ and $X^0$. Without loss of generality, we assume $X^0=0$.

For $n_0\in\N$, we take $a\ge\left(\P\left[A(n_0)\cdots A(1)\textrm{ has not rank~1 }\right]\right)^{1/n_0}$. But if $A(n)\cdots A(m)$ has not rank~1, then for any integer less than $\frac{n-m-n_0}{n_0}$, the operator $A(1+in_0)\cdots A((i+1)n_0)$ has not rank~1 either. From the independence of the $A(n)$, we deduce 
$$\P\left(A(n)\cdots A(m+1)\textrm{ has not rank~1 }\right)\le a^{n-m-n_0}.$$

We estimate $\delta\left(A(n)\cdots A(m+1)0,A(n)\cdots A(n_0+1+m)0\right)$: it is $0$ when $A(n+m)\cdots A(n_0+1+m)$ has rank~1, and it is always less than $\delta\left(A(n_0+m)\cdots A(m+1)0,0\right)$, that is less than $\I_{\{A(n)\cdots A(n_0+m+1)\textrm{ has not rank~1}\}} \left|A(n_0+m)\cdots A(m+1)0\right|_\mathcal{P}$, where $\I$ denotes the indicator function. Therefore, we have for any $n\ge m+n_0$
\begin{eqnarray}
\lefteqn{\E\left[\delta^l\left(A(n)\cdots A(m)0,A(n)\cdots A(n_0+1+m)0\right)\right]}\nonumber\\
&\le& \E\left[\I_{\{A(n)\cdots A(n_0+m+1)\textrm{ has not rank~1 }\}} \left|A(n_0+m)\cdots A(m+1)0\right|^l_\mathcal{P}\right]\nonumber\\
&=&a^{n-m-2n_0}\E\left[\left|A(n_0)\cdots A(1)0\right|^l_\mathcal{P}\right].
\end{eqnarray}
Let $n=qn_0+r$ be the Euclidean division of $n$ by $n_0$. Then we have
\begin{eqnarray*}
\left|x(n,0)\right|_\mathcal{P}&=&\delta\left(A(n)\cdots A(1)0,0\right)\\
&\le&\sum_{i=1}^{q}\delta\left(A(n)\cdots A(in_0+1)0,A(n)\cdots A((i-1)n_0+1)0\right)\\&&+\delta\left(A(n)\cdots A(n-r+1)0,0\right).\\
\end{eqnarray*} 
Therefore we have:
\begin{eqnarray}
\left(\E\left[\left|x(n,0)\right|^l_\mathcal{P}\right]\right)^{1/l}
&\le&\sum_{i=1}^{q}\left(a^{n-in_0-2n_0}\E\left[\left|A(n_0)\cdots A(1)0\right|^l_\mathcal{P}\right]\right)^{1/l}\nonumber\\
&&+\left(\E\left[\left|A(r)\cdots A(1)0\right|^l_\mathcal{P}\right]\right)^{1/l}.\label{eq10}
\end{eqnarray}
We apply this decomposition again (with $n=r$, $n_0=1$ and $a=1$), to check that 
$$ \left(\E\left[\left|A(r)\cdots A(1)0\right|^l_\mathcal{P}\right]\right)^{1/l}\le r \left(\E\left[\left|A(1)0\right|^l_\mathcal{P}\right]\right)^{1/l}\le n_0\left(\E\left[\left|A(1)0\right|^l_\mathcal{P}\right]\right)^{1/l} .$$
It follows from the MLP property, that there is $n_0\in \N$ such that $a<1$. Introducing the last equation in equation~(\ref{eq10}), we see that $$\left(\E\left[\left|x(n,0)\right|^l_\mathcal{P}\right]\right)^{1/l}\le\left(1+\frac{a^{-2n_0l}}{1-a^{n_0l}}\right)n_0\left(\E\left[\left|A(1)0\right|^l_\mathcal{P}\right]\right)^{1/l}.$$
\end{proof}

To go from the abstract LLT and renewal theorem to ours, we will use the following classical approximation lemma.
\begin{lem}\label{lemmsuppcomp}
Let $h$ be a continuous function with compact support from \mbox{$\R^a\times\R^b$.} Then there are two continuous functions $f_0$ and $g_0$ with compact support in $\R^a$ and $\R^b$ respectively, so that for any $\epsilon>0$, there are $f_i$ and $g_i$ continuous functions with compact support satisfying:
$$\forall x\in\R^a,y\in\R^b, |h(x,y)-\sum_if_i(x)g_i(y)|\le \epsilon f_0(x)g_0(y).$$
\end{lem}
In the sequel, we denote by $\mathcal{L}$ the Lebesgue measure.

\begin{proof}[Proof theorem~\ref{TLL}]
By theorem~\ref{TCL}, the algebraic non arithmeticity ensures that $\sigma>0$.

We apply proposition~\ref{Hm} and theorem~C of~\cite{HH1} or proposition~\ref{HS} and theorem~C of~\cite{HH2}. This proves that, if $g\in\mathcal{C}_c(\R)$ and if $f$ is a bounded Lipschitz function on $\prm$, then for any $x^0\in\R^d$,
\begin{equation}
\lim_n\sup_{u\in\R}\left|\sigma\sqrt{2\pi n}\E\left(f\left(\bar{x}(n,x^0)\right)g\left(\phi\left(x(n,x^0)\right)-\phi\left(x^0\right)-u\right)\right)-e^{-\frac{u^2}{2n\sigma^2}}\nu_0(f)\mathcal{L}(g)\right|=0.
\end{equation} 

Moreover these convergences are uniform in  $x^0$, because $\delta_{\overline{x^0}}$ is bounded independently of $x_0$ as a linear form on $\B^\infty$ and is in a disk of $\B_{\eta}$ with center $0$ and radius $\|X^0\|_\infty$ if $|x^0|_\mathcal{P}\le\|X^0\|_\infty$.
The uniformity allows us to take any random initial condition $X^0$ and get
\begin{equation}\label{eqTLL}
\lim_n \sup_{u\in\R}\left|\sigma\sqrt{2\pi n}\E\left[f\left(\bar{x}(n,X^0)\right)g\left(\phi\left(x(n,X^0)\right)-u\right)\right]-\E\left[e^{-\frac{(u+\phi(X^0))^2}{2n\sigma^2}}\right]\nu_0(f)\mathcal{L}(g)\right|=0.
\end{equation}

But the density of bounded Lipschitz functions in $(\mathcal{C}_c(\prm),\|.\|_\infty)$ allows us to take $f$ and $g$ continuous functions with compact support in equation~(\ref{eqTLL}). Now, it follows from lemma~\ref{lemmsuppcomp}, that for any $h$ continuous with compact support
\begin{equation}
\lim_n \sup_{u\in\R}\left|\sigma\sqrt{2\pi n}\E\left[h\left(\bar{x}(n,X^0),\phi\left(x(n,X^0)\right)-u\right)\right]-\E\left[e^{-\frac{(u+\phi(X^0))^2}{2n\sigma^2}}\right]\nu_0\otimes \mathcal{L}(h)\right|=0.
\end{equation}
According to lemma~\ref{bilip} the function $\Phi:x\mapsto (\phi,\bar{x})$ is Lipschitz with a Lipschitz inverse, therefore $h$ has  compact support iff $h\circ\Phi$ does. Since $\overline{x+u\1}=\overline{x}$, this concludes the proof.
\end{proof}

\begin{proof}[Proof of proposition~\ref{optNA}]
Assume the sequence of random variables is algebraically arithmetic and the conclusion of theorem~\ref{TLL} holds.

There are $a,b\in\R$ and $\theta$ with rank~1, such that every $A\in S_A$ and $\theta'\in T_A$ with rank~1 satisfy equation~(\ref{eqANA}). We set $t=\frac{2\pi}{b}$ if $b\neq 0$ and $t=1$ otherwise, and for any $x\in\R^d$, $w(\overline{x})=e^{it\left(\phi(\theta x)-\phi(x)\right)}$.
Equation~(\ref{eqANA}) implies that, for any $A\in S_A$, $y\in\R^d$, and any $\theta'\in T_A$ with rank~1:
$$e^{it\left(\phi(A \theta'y)-\phi(\theta'y)\right)}w(\overline{A\theta'y})=e^{ita}w(\overline{\theta'y}).$$
We chose $y$ such that $\phi(\theta'y)=0$.
By induction, we get
\begin{equation}\label{eqw}
e^{it\phi\left(x(n,\theta'y)\right)}w(\overline{x}(n,\theta'y))=e^{itna}w(\overline{\theta'y}).
\end{equation}

For any $f:\R\mapsto \R$ and $g:\prm\mapsto \R$ continuous with compact supports, the conclusion of theorem~\ref{TLL} for $h$ defined by $h(x)=f(\phi(x))(gw)(\overline{x})$ is that:
$$
\sigma\sqrt{2\pi n}\E\left[f(\phi(x(n,\theta'y)))g(\overline{x}(n,\theta'y))w(\overline{x}(n,\theta'y))\right]\rightarrow \mathcal{L}(f)\nu_0(gw). $$
Together with equation~(\ref{eqw}), it means
\begin{equation}\label{eq6}
e^{itna}w(\overline{\theta'y})\sigma\sqrt{2\pi n}\E\left[e^{-it.}f(\phi(x(n,\theta'y)))g(\overline{x}(n,\theta'y))\right]\rightarrow \mathcal{L}(f)\nu_0(gw)
\end{equation}

But conclusion of theorem~\ref{TLL} for $h$ defined by $h(x)=(fe^{-it.})(\phi(x))g(\overline{x})$ is that:
\begin{equation}\label{eq7}
\sigma\sqrt{2\pi n}\E\left[e^{-it.}f(\phi(x(n,\theta'y)))g(\overline{x}(n,\theta'y))\right]\rightarrow \mathcal{L}(fe^{-it.})\nu_0(g)
\end{equation}

Equations (\ref{eq6}) and (\ref{eq7}) together imply that $ta\in 2\pi\Z$ and that 
$$w(\overline{\theta'y})\mathcal{L}(fe^{it.})\nu_0(g)=\mathcal{L}(f)\nu_0(gw).$$
The right side of the equation does not depend on $\theta'$ so by lemma~\ref{suppnu} $w$ is constant on $S_{\nu_0}$, this proves $\nu_0(gw)=w(\overline{\theta'y})$, so $\mathcal{L}(fe^{-it.})=\mathcal{L}(f)$ that is $e^{it.}=1$ or $t=0$. This is a contradiction, which concludes the proof.
\end{proof}

\begin{proof}[Proof of theorem~\ref{renouv}]
Applying proposition~\ref{Hm} and theorem~D of~\cite{HH1}, we have that, if $g\in\mathcal{C}_c(\R)$ and if $f$ is a bounded Lipschitz function on $\prm$, then for any $x^0\in\R^d$,
$$\lim_{a\rightarrow-\infty} \sum_{n\ge1}\E\left[f\left(\phi\left(x(n,x^0)\right)-\phi\left(x^0\right)-a \right)g\left(\bar{x}(n,x^0)\right)\right]=0,$$
$$ \lim_{a\rightarrow+\infty} \sum_{n\ge1}\E\left[f\left(\phi\left(x(n,x^0)\right)-\phi\left(x^0\right)-a \right)g\left(\bar{x}(n,x^0)\right)\right]=\frac{\nu_0(f)\mathcal{L}(g)}{\gamma}.$$
Moreover these convergences are uniform in  $x^0$, because $\delta_{\overline{x^0}}$ is bounded as a linear form on $\B^\infty$.
The uniformity allows us to remove the $\phi\left(x^0\right)$ in the last equations and take any random initial condition.
The result follows by the same successive approximations as in the proof of the LLT.
\end{proof}

\begin{proof}[Proof of theorem~\ref{PGD}]
Without lost of generality, we can assume that $\gamma=0$.

The exponential moment of $\sigma(A)$ means that there is a $\theta>0$ such that \mbox{$\int e^{\theta \sigma(A)} d\mu(A)<\infty$.}
An easy bound of the norm of $\xi^k(y)Q(.,dy)$ inspired by the proof of proposition~\ref{Hm} ensures that $z\mapsto Q_z$ is analytic on the open ball with center $0$ and radius $\theta$.
To prove that it is continuous on the domain $\{|\mathcal{R}z|<\theta/2\}$, we apply the same method.

Now theorem~E of~\cite{HH1} gives $$\lim_{n}\frac{1}{n}\ln \P\left[\phi\left(x(n,X^0)\right)-\phi(X^0)>n\epsilon\right]=-c(\epsilon).$$
Let $0<\eta<\epsilon$. For any $n\ge \|\phi(X^0)/\eta\|_\infty$, we have
$$\P\left[\phi\left(x(n,X^0)\right)-\phi(X^0)>n\epsilon\right]
\ge \P\left[\phi\left(x(n,X^0)\right)>n(\epsilon+\eta)\right],$$
which implies that
$$\liminf_{n}\frac{1}{n}\ln \P\left[\phi\left(x(n,X^0)\right)-n\gamma>n\epsilon\right]\ge-c(\epsilon+\eta).$$
The same method gives
$$\limsup_{n}\frac{1}{n}\ln \P\left[\phi\left(x(n,X^0)\right)-n\gamma>n\epsilon\right]\le-c(\epsilon-\eta).$$
By continuity of $c$, the first equality is proved. The second one follows from the same method applied to $-\phi$ instead of $\phi$.
\end{proof}
\subsection{Max-plus case}
Before proving the statements, we recall a few needed definitions and results about powers of matrices in the \mp algebra.
\begin{defn}
\begin{enumerate}
\item The critical graph of $A$ is obtained from $\mathcal{G}(A)$ by keeping only nodes and arcs belonging to circuits with average weight $\rm(A)$. It will be denoted by $\mathcal{G}^c(A)$.
\item The cyclicity of a graph is the greatest common divisor of the length of its circuits if it is strongly connected (that is if any node can be reached from any other). Otherwise it is the least common multiple  of the cyclicities of its strongly connected components. The cyclicity of $A$ is that of $\mathcal{G}^c(A)$ and is denoted by $c(A)$.
\end{enumerate}
\end{defn}

\begin{rem}\label{powerint} Interpretation of powers with $\mathcal{G}(A)$.\\
If $(i_1,i_2\cdots,i_n)$ is a path on $\mathcal{G}(A)$, its weight is $\sum_{1\le j\le n-1}A_{i_ji_{j+1}}$, so that $\left(A^{ n}\right)_{ij}$ is the maximum of the weights of length $n$ paths from $i$ to $j$.
\end{rem}

\begin{thm}[\cite{cohen83}]\label{proppuiss}
Assume $\mathcal{G}(A)$ is strongly connected, \mbox{$\rm(A)=0$}. Then the sequence $\left(A^{ n}\right)_{n\in\N}$ is ultimately periodic and the ultimate period is the cyclicity of $A$.
\end{thm}

\begin{proof}[Proof of theorem~\ref{s>0}]
Suppose that $\sigma=0$. By proposition~\ref{Hm} we may apply theorem~A of~\cite{HH1}. The third point of the theorem says that there exists a bounded Lipschitz function $f$ such that for $\nu$-almost every $(B,\bar{x})$:
\begin{equation}\label{cobord}
\max_i(Bx)_i-\max_ix_i=f(\overline{x})-f(\overline{Bx})
\end{equation}
Since all functions in that equation are continuous, every $B\in S_A$ and $\overline{x}\in Supp(\nu_0)$ satisfy equation~(\ref{cobord}). If $B\in S_A$ and $\overline{x} \in Supp(\nu_0)$, then $\overline{Bx}\in Supp(\nu_0)$, so by induction equation~(\ref{cobord}) is satisfied by any $B$ in $T_A$. Since for $B\in T_A$, $B^n\in T_A$, $\max_iB^nx_i$ is bounded. But there exists a $k$ such that $c(B)\rm(B)=B^{c(B)}_{kk}$, so $\max_i\left(B^{nc(B)}x\right)_i\ge nc(B)\rm(B) +x_k$ and $\rm(B)\le 0$.\\

Since every path on $\mathcal{G}(B)$ can be split into a path with length at most $d$ and closed paths whose average length are at most $\rm(B)$, we have:
$$\max_i(B^nx)_i\le (n-d)\rm(B) +d \max_{B_{ij}>-\infty} |B_{ij}|+\max_ix_i,$$
therefore $\rm(B)\ge 0$.

So $\sigma=0$ implies that $\forall B\in T_A,\rm(B)=0$.\\

Conversely, if $\rm(B)=0$ for every $B\in T_A$, then 
\begin{eqnarray*}
\max_ix_i(n,0)&=&\max_{ij}\left(A(n)\cdots A(1)\right)_{ij}\\
&\ge& \rm\left(A(n)\cdots A(1)\right)=0 \textrm{ a.s.\,}.
\end{eqnarray*}
Therefore $\mathcal{N}(0,\sigma^2)(\R_+)\ge1$, and $\sigma=0$. 
\end{proof}

\begin{proof}[Proof of theorem~\ref{NA}]
We assume the system is algebraically arithmetic. 
Then there are $a,b\in\R$ and $\theta\in T_A$ such that for any $A \in S_A$ and $\theta'\in T_A$ with rank~1, we have:
$(\theta A\theta'-\theta \theta')(\R^d) \subset (a+b\Z)\1.$
Replacing $\theta'$ by $A^n\theta'$, we get $(\theta A^{n+1}\theta'-\theta A^{n} \theta')(\R^d) \subset (a+b\Z)\1$ and by induction 
\begin{equation}\label{eq4}
(\theta A^{n+k}\theta'-\theta A^{n} \theta')(\R^d) \subset (ka+b\Z)\1
\end{equation}
From now on, we assume that $\mathcal{G}(A)$ is strongly connected.
The matrix $\tilde{A}$ defined by $\tilde{A}_{ij}=A_{ij}-\rm(A)$ satisfy $\rm\left(\tilde{A}\right)=0$ and has a strongly connected graph. Therefore, by theorem~\ref{proppuiss}, there is an $n$ such that for any  indices~$i,j$, $\tilde{A}^{n+c(A)}_{ij}=\tilde{A}^n_{ij}$.

Since for any $n\in\N$, $A^n_{ij}=\tilde{A}^n_{ij}+n\rm(A)$, it means that $A^{(n+1)c(A)}_{ij}=A^{nc(A)}_{ij}+c(A)\rm(A)$, and  $(\theta A^{n+c(A)}\theta')_{ij}-(\theta A^{n} \theta')_{ij}=c(A)\rm(A)$.
Together with equation~(\ref{eq4}), it says that $\rm(A)\in a+\frac{b}{c(A)}\Z\subset a+\frac{b}{d!}\Z$, which concludes the proof.
\end{proof}
\section{Acknowledgements}
The author gratefully thanks Jean Mairesse for useful talks and suggestions of improvements to this article.
\bibliographystyle{alpha}
\bibliography{max+}
\end{document}